\documentclass[journal]{IEEEtran}
\IEEEoverridecommandlockouts
\usepackage{balance} 
\balance
\usepackage{graphics} 
\usepackage{graphicx}
\usepackage{epsfig} 
\usepackage{times} 
\usepackage[tbtags]{amsmath}
\usepackage{amssymb}  
\usepackage{stfloats}
\usepackage{array}
\usepackage{bm}
\DeclareMathOperator*{\Min}{Min}
\usepackage{float}
\usepackage{booktabs, tabularx, multirow, array}
\usepackage{nomencl}
\makenomenclature
\usepackage{siunitx}
\usepackage{MnSymbol,wasysym}
\allowdisplaybreaks
\DeclareSIUnit\kWh{kWh}
\DeclareSIUnit\DolarMWh{\$ /MWh}
\usepackage{calligra}
\DeclareMathAlphabet{\mathcalligra}{T1}{calligra}{m}{n}

\usepackage{subcaption}
\def\x{\mathbf{x}}
\def\y{\mathbf{y}}
\def\d{\mathbf{d}}
\def\D{\mathbf{D}}
\def\L{\mathbf{L}}
\def\e{\mathbf{e}}

\usepackage{dcolumn,booktabs}
\newcolumntype{d}[1]{D{.}{.}{#1}}

{\LARGE \bf
\title{Stochastic Time-of-use-type Constraints for Uninterruptible Services}
}

\author{Ana~Batista,
        David~Pozo,~\IEEEmembership{Senior Member,~IEEE,}
        and~Jorge~Vera 
\thanks{This work was supported by Skoltech NGP Program (Skoltech-MIT joint project). A. Batista is with the Department of Industrial Engineering,
        Pontificia Universidad Cat\'olica de Chile, 
        Santiago, Chile and Skolkovo Institute of Science and Technology, Moscow, Russia (email: a.batista@skoltech.ru). D. Pozo  is with the Center for Energy Systems, Skolkovo Institute of Science and Technology (Skoltech), Moscow, Russia (email: d.pozo@skoltech.ru). J. Vera  is with the Department of Industrial Engineering,  Pontificia Universidad Cat\'olica de Chile, 
        Santiago, Chile (email: jvera@uc.cl)}   
        }
\begin{document}

\maketitle

\begin{abstract}
In this paper, a mixed integer linear formulation for problems considering time-of-use-type constraints for uninterruptible services is presented. Our work is motivated by demand response problems in power systems, in which certain devices must remain online once they are switched on, along with a certain number of periods. Classically, this kind of constraints are modeled as a summation over a rolling time windows. This makes it difficult to consider this time-of-use parameter  as uncertain. We propose an alternative formulation in which the time of use is on the right-hand side of a constraint instead on the limit of a summation. This allows applying existing stochastic optimization methodologies easily. An illustrative model for the optimal allocation of an uninterruptible load for the demand response problem supports the proposed formulation.
\end{abstract}


\section{Introduction}
Many optimal allocation problems incorporate decisions about service (e.g., power supply) duration that once it is started cannot be interrupted along with a determinate time window (minimum time of use). There is also the possibility that the uninterruptible service has an expiration time (maximum time for using) or a determined time window (exact time of use). There exists a vast amount of research in several fields that include formulations for the \mbox{\textit{time-of-use}} of uninterruptible services by means of integer variables.
For instance, the unit commitment problem seeks to allocate power units considering integer minimum up-time and down-time constraints. It is classically reformulated as mixed linear constraints using summation of binary variables over a time windows \cite{carrion2006computationally,ostrowski2012tight,ozturk2004solution}. A similar framework can be found in other areas, such as healthcare. For example, in reference \cite{bachouch2012integer}, the patient allocation is constrained by an integer length of stay, over a time windows. 
Demand response (DR) problems, \cite{deng2015survey}, \cite{albadi2008summary} also include \mbox{{time-of-use-type}} constraints. This problem intends to shape the user consumption by incentivizing load shifting across hours. In brief, the user will prefer to use a flexible load while the electricity price is low. Therefore, different loads should be allocated along the time horizon, e.g., a day. 
For instance, \cite{chen2012real} developed a MILP real-time price-based DR approach applied to home appliances. They studied, in particular, uninterruptible and deferrable home appliances in which the time-of-use parameter is deterministic. The time-of-use is modeled as a deterministic time window over the indexes of a summation, making it difficult to consider uncertainty in the operation time duration.

In general, previous works from the existing literature have considered time-of-use-type constraints as an integer parameter, while the nature of this parameter is stochastic in many applications. In this letter, we propose a new simple, but effective formulation that includes constraints for problems in which a service cannot be interrupted once it has been started. The service has an uncertain time to be fulfilled. The resulting model has the time-of-use parameter on the right-hand side of the constraints, rather than over the indexes of a summation. This structure facilitates the implementation of the existing algorithms (e.g., dual-based methods, or bender decomposition) that considers uncertainty, such as stochastic programming, robust optimization, and distributionally robust optimization.

\section{Constraints Formulation}
\label{sec:3}
\begin{figure} [h]
\centering
\includegraphics[scale=0.55]{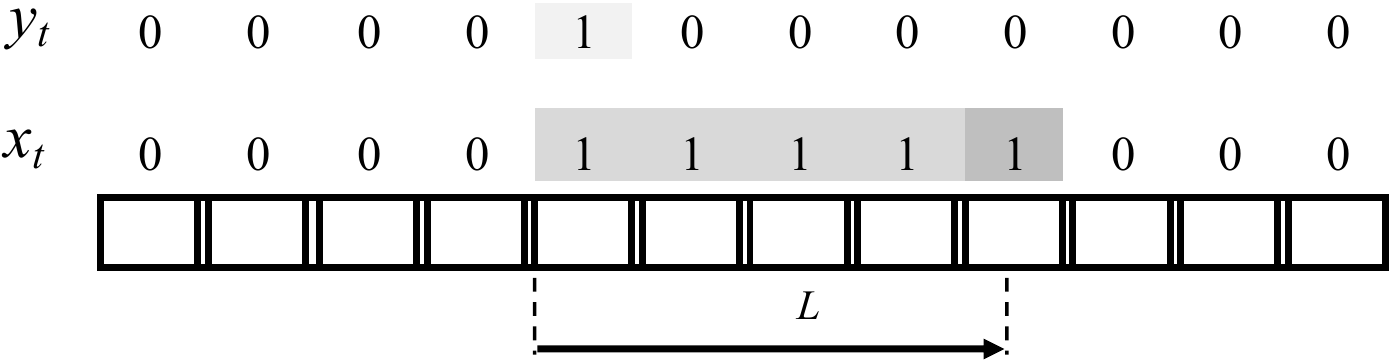}
\caption{{An illustrative representation of the allocation problem for uninterruptible service}}
\label{fig:1}
\end{figure}

The proposed framework is sketched in Figure \ref{fig:1}. The allocation process is subdivided into a set of $T$ time slots. The variable $y_t$ indicates the initiation of the service. The binary variable $x_t$ associated with each time slot takes the value 1 if we allocate this slot to fulfill the service, and 0 otherwise. The parameter $L$ is the continuous time of use of the uninterruptible service. Without loss of generality, we assume that the service always initiates at the beginning of a time slot and it is assigned completely; thus it can end at any moment and not necessarily at the end of a time slot. The service has to be allocated only once during the time horizon, and there must be enough time slots to complete the service length $L$.

The proposed set of linear constraints that represent the aforementioned  problem description can be formulated as follows:

\newpage
\begin{IEEEeqnarray}{l"l}
	\sum_{t=1}^{T}  y_t = 1  \label{eq:1} \\  
	\sum_{t=1}^{T}  x_t \geq L  \label{eq:2} \\ 
	y_t \geq x_t - x_{t-1}   &   \qquad  \forall t=\{2,\ldots, T\} \label{eq:3}  \\ 
	y_t \geq x_t    &   \qquad  \forall t=1  \label{eq:4}\\
	x_t \in \{0,1\},\ y_t \in [0,1] &   \qquad \forall  t\in T \label{eq:5}
\end{IEEEeqnarray}

Constraint \eqref{eq:1} indicates that the service has to be assigned once during the entire horizon. 
Constraint \eqref{eq:2} guarantees that the service has enough time slots allocated to accommodate the service of length $L$. Constraint \eqref{eq:3} depicts the logic between the allocation variable $x_t$ and the activation variable $y_t$. 
We assume that for $t=0$, $x_0 = 0$. Thus, for $t=1$, $y_t \geq x_t$, as defined in constraint \eqref{eq:4}. 
It should be noted that constraints \eqref{eq:1} and \eqref{eq:2}, ensure that the service cannot be interrupted and there are enough time slots to be allocated. 
Lastly, constraint \eqref{eq:5}, defines the variables' domain. Observe that it is not required to define $y_t$ as binary due to equations \eqref{eq:1}, \eqref{eq:3}, and \eqref{eq:4} enforce to $y_t$ to take the value $1$ for a single time slot and therefore it would be $0$ for the rest of time slots. 

\section{Optimal Load Allocation Model for DR}
\subsection{Deterministic Model}
\label{sec:A}
In this section, we present an illustrative model for the optimal load allocation for the DR problem. The time-of-use-type constraints are related to the time required for fulfilling a load (e.g., electric vehicle charging when the initial state of charge of the batteries is not known). 
For simplicity, we consider only one load to be allocated along 24 hours. We assumed time slots of 1 hour indexed by $t \in [1, 24]$. The time-of-use is assumed to be non-integer and fixed over the time horizon $T=24$.
We consider that the maximum and minimum power is constant during all interval $L$. The maximum power of the load (e.g., maximum battery charging power) is given by $\overline{r}$, while the minimum power of the load is restricted to $\underline{r}$ for every time slot that the load demand is scheduled. 
Thus, we define $d_t$ as the energy allocated for each time slot, and $\lambda_t$ as the cost per unit of energy for each hour (typically, the price of electricity). The mathematical formulation for this problem is described in \eqref{eq:6}--\eqref{eq:9}.

\begin{align}
&\Min_{d,x,y}
\begin{aligned}[t]
&\displaystyle\sum_{t=1}^{T}  \; \lambda_{t} d_{t}  \label{eq:6}
\end{aligned}  \\ 
&\text{s.t.:}  \notag \\ 
& x_t \underline{r} \leq d_t \leq x_t \overline{r} &\forall t \in T \label{eq:7}\\ 
&\displaystyle\sum_{t=1}^{T} d_t = D  \label{eq:8} \\
&\displaystyle\sum_{t=1}^{T}  y_t = 1  \label{eq:1.1} \\  
&\displaystyle\sum_{t=1}^{T}  x_t \geq L  \label{eq:2.1} \\ 
& y_t \geq x_t - x_{t-1}   &   \qquad  \forall t=\{2,\ldots, T\} \label{eq:3.1}  \\ 
& y_t \geq x_t    &   \qquad  \forall t=1  \label{eq:4.1}\\
&	x_t \in \{0,1\},\ y_t \in [0,1] &\forall    t\in T \label{eq:5.1} \\
& d_t \geq 0    &\forall   t\in T  \label{eq:9}
\end{align}

The objective \eqref{eq:6} minimizes the total load allocation cost, which is defined as the sum of the allocated demand over the time period $t$. Constraint \eqref{eq:7} ensure that the hourly load demand is within its minimum and maximum limits. Constraint \eqref{eq:8} ensures the allocation of the total demand of the load for all periods. The group of equations \eqref{eq:1.1}--\eqref{eq:5.1} are the time-of-use-related constraints described in Section \ref{sec:3}.
Constraint \eqref{eq:9} sets the bound to the variable $d_t$. 

\subsubsection*{Numerical Results} 
\begin{figure}[t]
  \begin{subfigure}[t]{.495\linewidth}
    \centering
    \includegraphics[width=\linewidth]{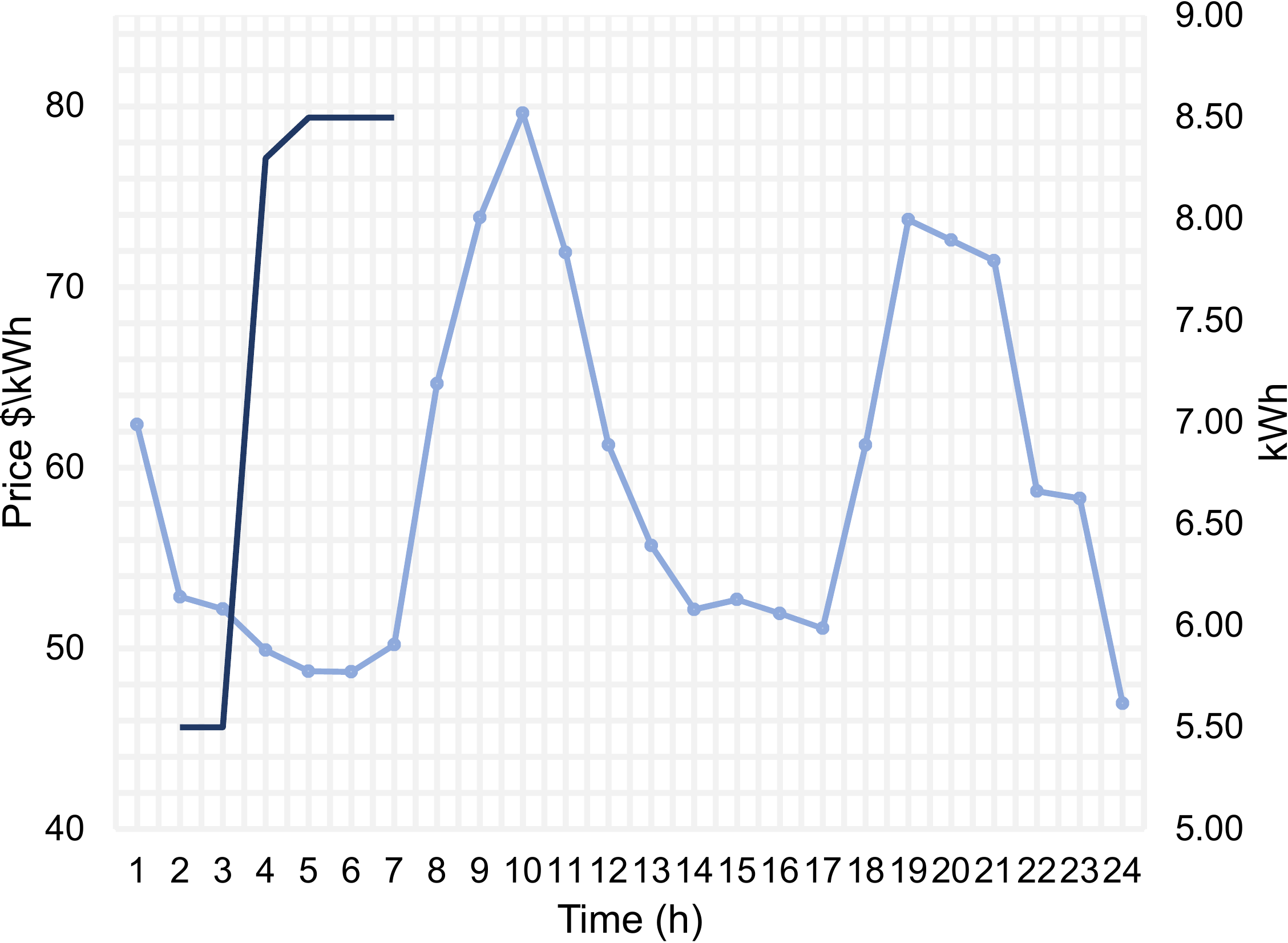}
    \caption{Case 1, L= 5.27 h}
  \end{subfigure}
  \hfill
  \begin{subfigure}[t]{.495\linewidth}
    \centering
    \includegraphics[width=\linewidth]{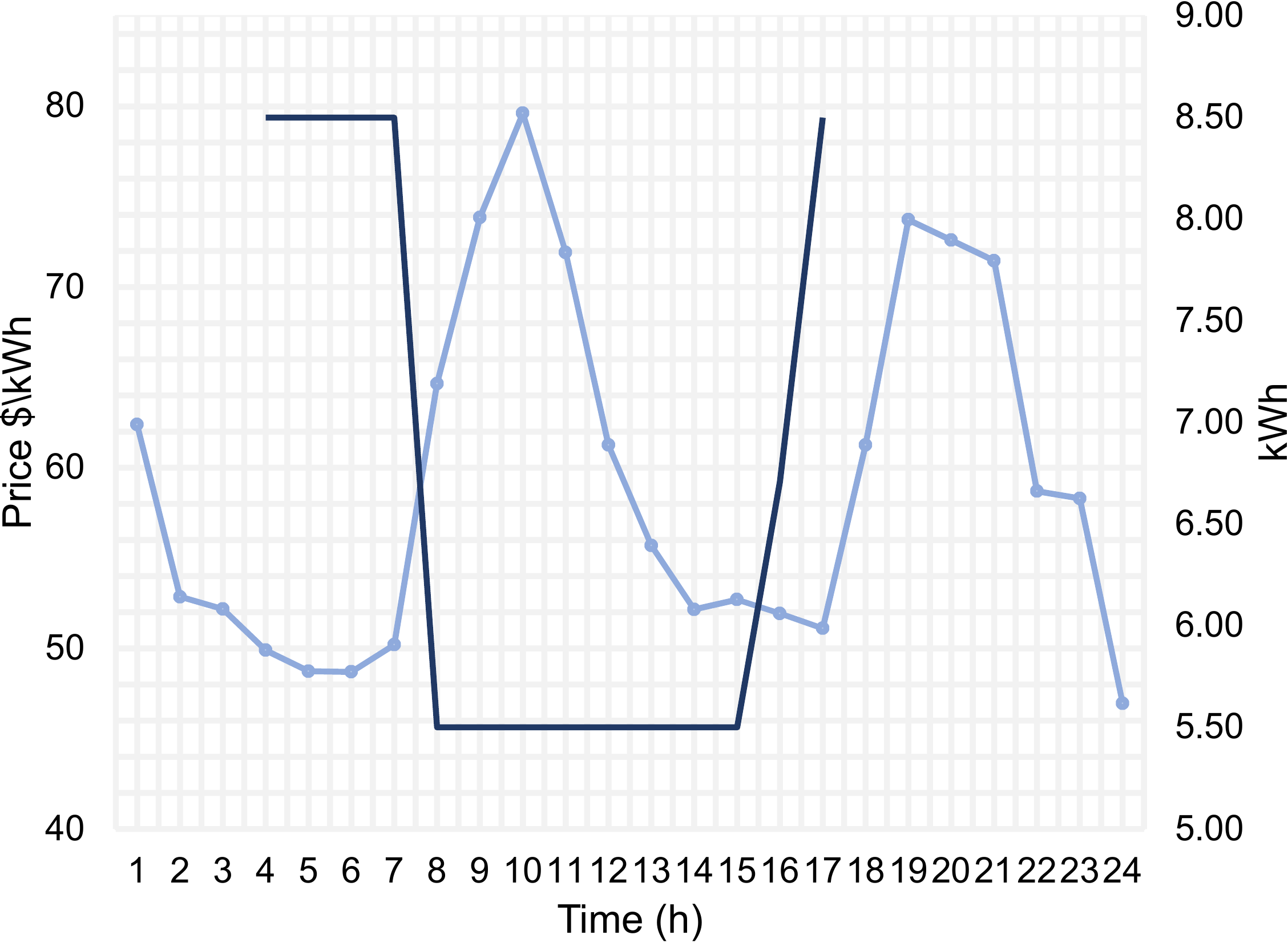}
    \caption{Case 2, L= 10.97 h}
  \end{subfigure}

  \medskip

  \begin{subfigure}[t]{.495\linewidth}
    \centering
    \includegraphics[width=\linewidth]{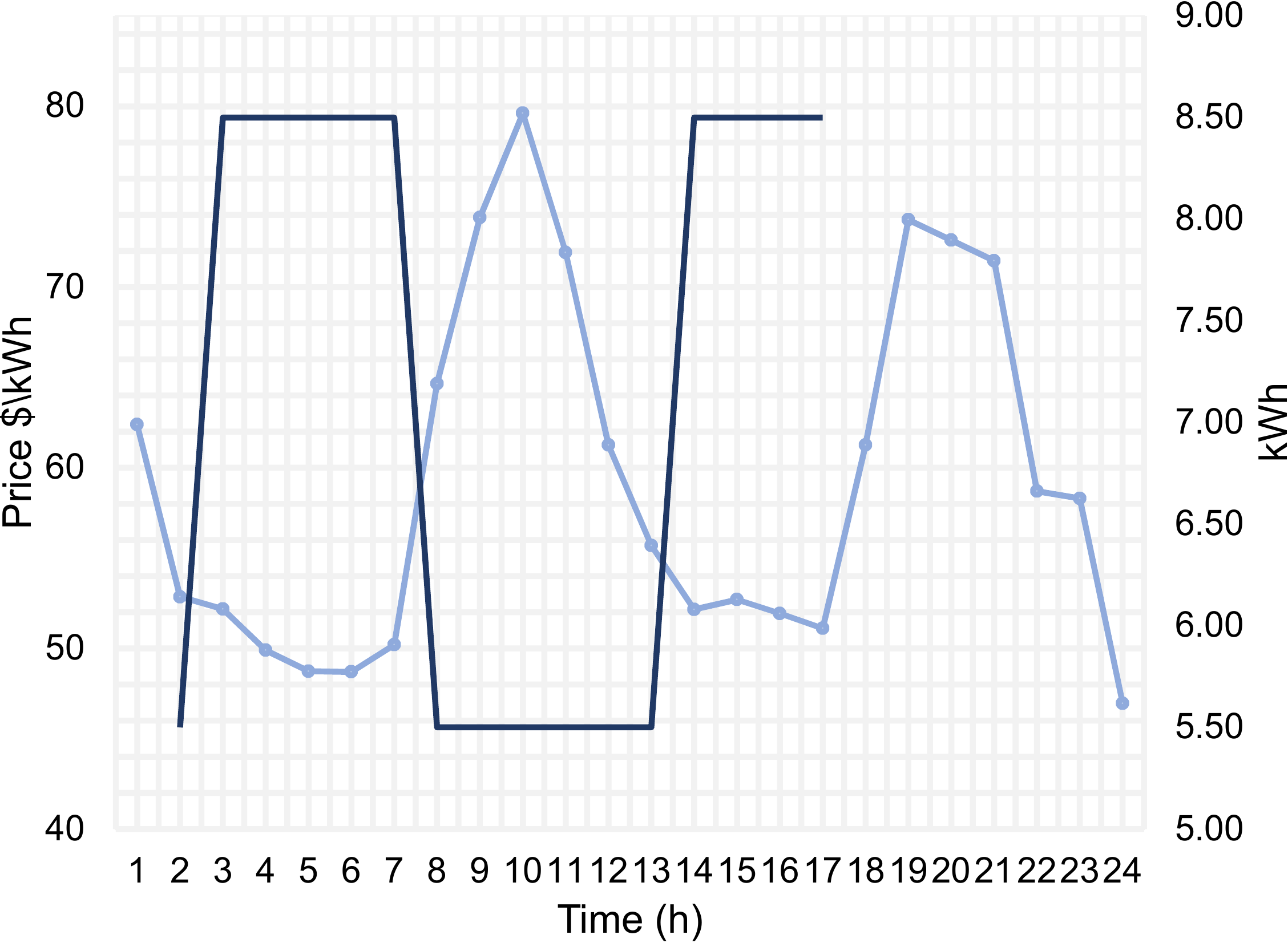}
    \caption{Case 3, L= 13.53 h}
  \end{subfigure}
  \hfill
  \begin{subfigure}[t]{.495\linewidth}
    \centering
    \includegraphics[width=\linewidth]{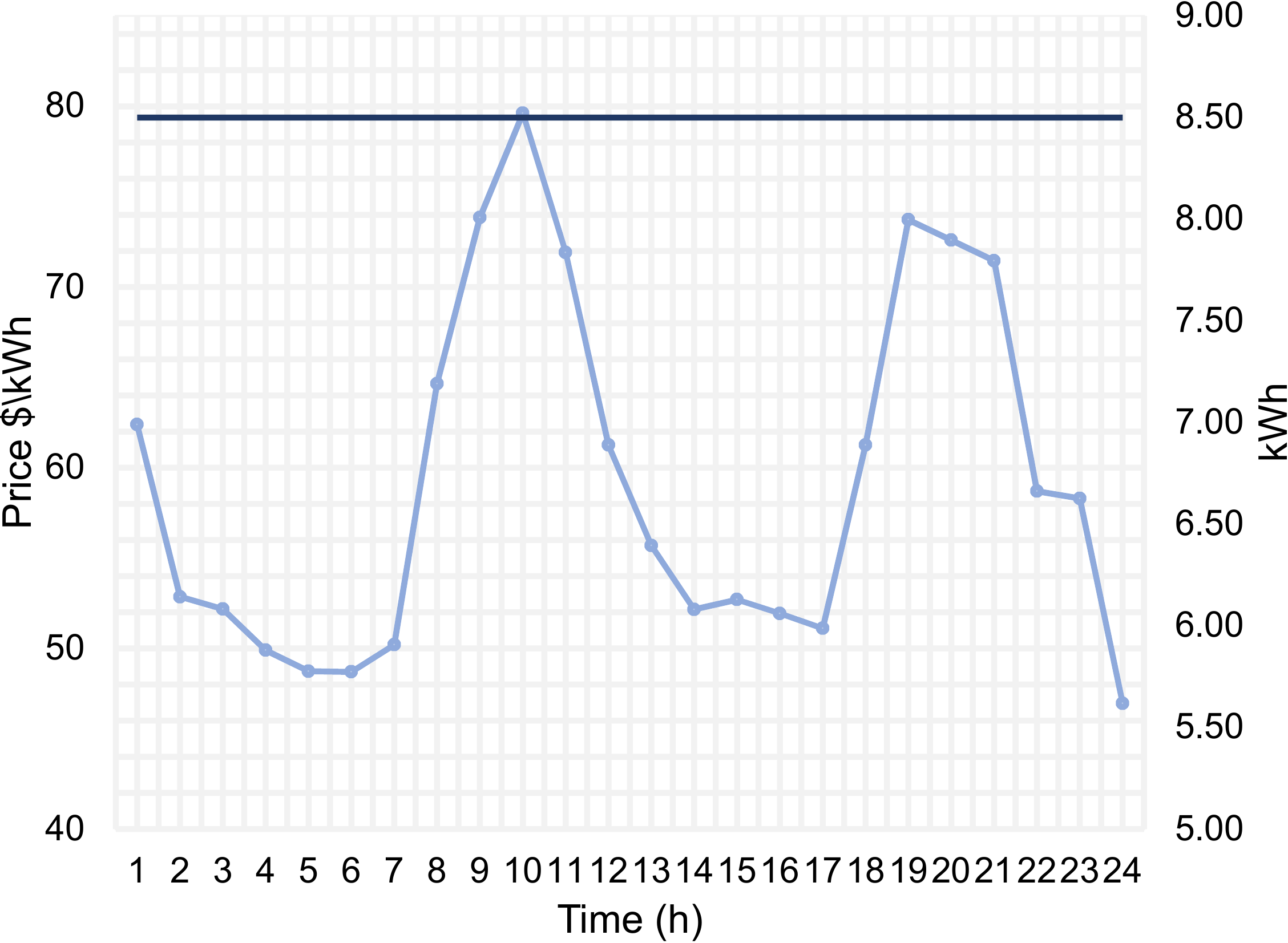}
    \caption{Case 4, L= 23.33 h}
  \end{subfigure}
  \caption{Results of the hourly allocation of the load compared to the price curve}
\label{fig:2}
\end{figure}

We employed price data of the load to be allocated on an hourly basis. For a time horizon of 24 hours, we assumed two peaks of price from 8:00 to 12:00 hours and from 18:00 to 22:00. We considered a maximum constant ratio of load, $\overline{r}$, of 8.5 \si{\kW} and a minimum ratio of load, $\underline{r}$ of 5.5 \si{\kW}. The data employed for the demand and time-of-use is shown in Table \ref{table:1}. 
The results reported in Table \ref{table:1} and Figure \ref{fig:2} indicate that the model guarantees the allocation of the unit to the periods with the lower prices for all the cases. Considering that the time-of-use values are non-integers, the allocation reserves enough periods to accomplish the required load of charge at a minimum cost per period. Although the energy allocated varies along the hours (higher at lower prices), the load is not interrupted until it is satisfied.

\begingroup
\renewcommand{\arraystretch}{1.2} 
\begin{table}[htbp]
\centering
\caption{Data and optimal cost results}
\label{table:1}
\begin{tabular}{l|d{4.2}d{4.2}d{4.2}d{4.2}}
\hline
\multicolumn{1}{c|}{} &\multicolumn{1}{c}{C1} & \multicolumn{1}{c}{C2} &\multicolumn{1}{c}{C3} & \multicolumn{1}{c}{C4} \\ \hline
D (\si{\kWh}) & 44.80 & 93.25 & 115.00 & 198.30 \\ \hline
L (\si{h}) & 5.27 & 10.97 & 13.53 & 23.33 \\ \hline
Total Cost (\$) & 2\,250.50 & 5\,286.92 & 6\,428.50 & 11\,664.20 \\ \hline
\end{tabular}
\end{table}
\endgroup

\subsection{Robust Model}
Motivated by the fact that demand loads, and therefore, time-of-use parameters can be uncertain, we reformulated the deterministic model \eqref{eq:6}--\eqref{eq:9} as a robust model \eqref{eq.objRO}--\eqref{eq.defDL}. Also, we extended the model for more than one load. We consider uncertainty for total energy demand (e.g., the state-of-charge of an electric vehicle) represented by the vector $\mathbf{D}$ indexed by $i \in I$.  Therefore, the time-of-use (time for charging in the electric-vehicle example) is also uncertain and represented in here by the vector $\mathbf{L}$. Uncertainty time-of-use is related to total demand thought the minimum and maximum rates for each load. Thus, uncertainty can be defined over $\mathbf{D}$, only. In our case, we considered ellipsoidal uncertainty to incorporate correlation between loads.  
The uncertainty set $\mathrm{U}$ is defined in \eqref{eq:uncertainty}. 
The parameter $\widehat{\mathbf{D}}$ is the vector of the nominal (average) demand value for the $I$ loads, and $\mathbf{D}$ is the uncertain demand vector. The parameter $\mathbf{A}$ is a lower triangular matrix which is obtained by the Cholesky decomposition of the estimated covariance matrix $\mathbf{\Sigma}$ 
\footnote{Observe that the ellipsoidal uncertainty set \eqref{eq:uncertainty} is derived from its classical definition \cite{ben1999robust}, 
\mbox{$\mathrm{U} = \big\{\mathbf{D} \vert (\mathbf{D} - \widehat{\mathbf{D}})^T \mathbf{{\Sigma}}^{-1}(\mathbf{D} - \widehat{\mathbf{D}})^T \leq \beta^2 \big\}$}. By replacing \mbox{$\mathbf{{\Sigma}}^{-1} = (\mathbf{AA}^T )^{-1}$} in $\mathrm{U}$, we have \mbox{$\mathrm{U} = \big\{\mathbf{D} \vert  \: \lVert (\mathbf{D} - \widehat{\mathbf{D}})^T  \mathbf{A}^{-1} \lVert_1 \leq \beta \big\}$}. Now it is easy rewrite $\mathbf{D}$ as  \mbox{$\mathbf{D} = \widehat{\mathbf{D}} + \mathbf{A} \mathbf{e}$} where $\mathbf{e}$ is the component-wise vector of errors limited by $-1 \leq e_i \leq 1$, and conservatives parameter $\beta$ is fixed to 1.}. 
Note that the uncertainty definition by \eqref{eq:uncertainty} is a linear conservative approach of the classic quadratic representation that preserves the original correlation \cite{ben1999robust}.

\begin{align} 
&\mathrm{U} = \big\{ \mathbf{D} \in \mathbb{R}^I, \mathbf{e} \in \mathbb{R}^I \vert \: \mathbf{D}  = \widehat{\mathbf{D}} + \mathbf{A} \mathbf{e}, \:  \lVert \mathbf{e} \lVert_{\infty} \leq 1 \big \}   \label{eq:uncertainty} 
\end{align}

The resulting robust DR allocation problem should guarantee minimal cost while ensuring feasibility for any scenario within the uncertainty set $\mathrm{U}$. Thus, we can state the robust problem as a bi-level one described as follows:

\begin{IEEEeqnarray}{lr}
    \max \limits_{{\D,\L,\e}} \quad \min \limits_{{\x,\y,\d}} \displaystyle\sum_{i=1}^{I} \displaystyle\sum_{t=1}^{T} \lambda_{t} d_{it}   \label{eq.objRO}  \\
\qquad \quad \text{s.t.:} \: x_{it} \underline{r}_i \leq d_{it} \leq x_{it} \overline{r}_i & \forall i \in I, t \in T  \label{eq.iniInner} \\
\qquad \qquad \displaystyle\sum_{t=1}^{T} d_{it} = D_i  & \forall i \in I \\
\qquad \qquad \displaystyle\sum_{t=1}^{T} y_{it} = 1  & \forall i \in I \\
\qquad \qquad  y_{it} \geq x_{it} &  \forall i \in I, t=1\\
\qquad \qquad y_{it} \geq x_{it} - x_{it-1} & \forall i \in I, t=\{2,\ldots, T\} \\
\qquad \qquad \displaystyle\sum_{t=1}^{T}  x_{it} \geq L_i & \forall i \in I\\
\qquad \qquad x_{it} \in \{0,1\}  & \forall  i \in I,  t\in T \\
\qquad \qquad  y_t \in [0,1] \geq 0  & \forall  i \in I,  t\in T \vspace{2\jot} \\
\qquad \qquad d_{it} \geq 0    & \forall i \in I, t \in T 
\label{eq.finInner}\\
\text{s.t.:} \: D_i = \widehat D_i + \displaystyle\sum_{i' \in I} A_{i,i'} e_{i'}  &  \forall i \in I  \label{eq.D_def}   \\
\quad \frac{D_i}{\underline{r}_i} \leq L_i \leq \frac{D_i}{\overline{r}_i}  &  \forall i \in I  \label{eq.L_def} \\
\quad  -1 \leq e_i \leq 1  & \forall i \in I  \label{eq. defE}  \\
\quad  D_i, L_i \geq 0  & \forall i \in I   \label{eq.defDL}
\end{IEEEeqnarray}

The robust formulation is a bi-level max-min problem. The inner (lower-level) problem \eqref{eq.objRO}--\eqref{eq.finInner} represents the optimal allocation of loads for a given $\L$ and $\D$. 
This problem is indeed an extension of the deterministic one, \eqref{eq:6}--\eqref{eq:9}, for multiple loads. Decision variables of the inner problem are $\x,\y$ and $\d$. $\L$ and $\D$ are parameters for this problem.
The outer (upper-level) problem \eqref{eq.objRO}, \eqref{eq.D_def}--\eqref{eq.defDL} selects an adversarial scenario of time-of-use and demands within the uncertainty set definition. This selection infers the largest cost to the inner problem considering that it would optimally minimize the total load cost allocation. Decision variables of the outer problem are $\D, \L$, and $\e$. 
The objective function \eqref{eq.objRO} represents the worst-case optimal load allocation cost. Equation \eqref{eq.D_def} is the total uncertain demand for each load defined accordingly to the uncertainty set \eqref{eq:uncertainty}. Equation \eqref{eq.L_def} defines the limits of the time-of-use-type constraint. Thus, the uninterruptible time in hours that the load requires.  Equations \eqref{eq. defE} and \eqref{eq.defDL} defines the domain limits of the decision variables.

The robust max-min problem \eqref{eq.objRO}--\eqref{eq.defDL} contains  a MILP inner problem. The inner binary variable, $\x$, complicates merging outer and inner problems into a single one by dualization of the inner problem. To address this issue, it is possible to construct a column-and-constraint generation algorithm \cite{ZhaoZeng}. Similarly to \cite{ZhaoZeng}, we implemented an iterative algorithm for solving the robust DR problem ensuring a global optimal solution.

\subsubsection*{Numerical Results} 
The data about the prices remains the same as described for the deterministic model. We consider two loads. 
The range of the power of both loads, $\underline{r}_i$ and $\overline{r}_i$ are set to $[3.5, 4.5]$ and $[2.5,3.5]$. 
We assumed that the demand follows a multivariate  distribution with known mean and covariance matrix, $\widehat{\mathbf{D}}$, and $\mathbf{\Sigma}$. We set the nominal mean value of $\widehat D_1 = 30 \ \si{kWh}$ and $\widehat D_2 = 20 \ \si{kWh}$. The values of the covariance matrix $\mathbf{\Sigma}$ in the main diagonal are defined as, $\sigma^2_1 = 10^2 \ \text{and} \ \sigma^2_2 = 15^2$, and the off-diagonal entries are set to  $\rho \times \sigma_1 \sigma_2$. Where $\rho$ is a scalar representing the correlation between loads.

Table \ref{table:2} shows the optimal results for values of correlation ranging from $\rho = [-1, +1]$ in which $\rho = 0$ corresponds to the uncorrelated case.
Observe that the maximum cost (worst-case demand, and therefore time-of-use) is for the case of no-correlation. In this case, the outer problem has a trivial solution (maximum load demand for both cases). When the correlation is considered, the solution is not trivial anymore. Note that the problem can take advantage of this information that would help to decrease the overall worst-case total cost about the no-correlation case, no matter if the correlation is positive or negative.

\begin{figure} [t]
\centering
\includegraphics[width=0.9\linewidth]{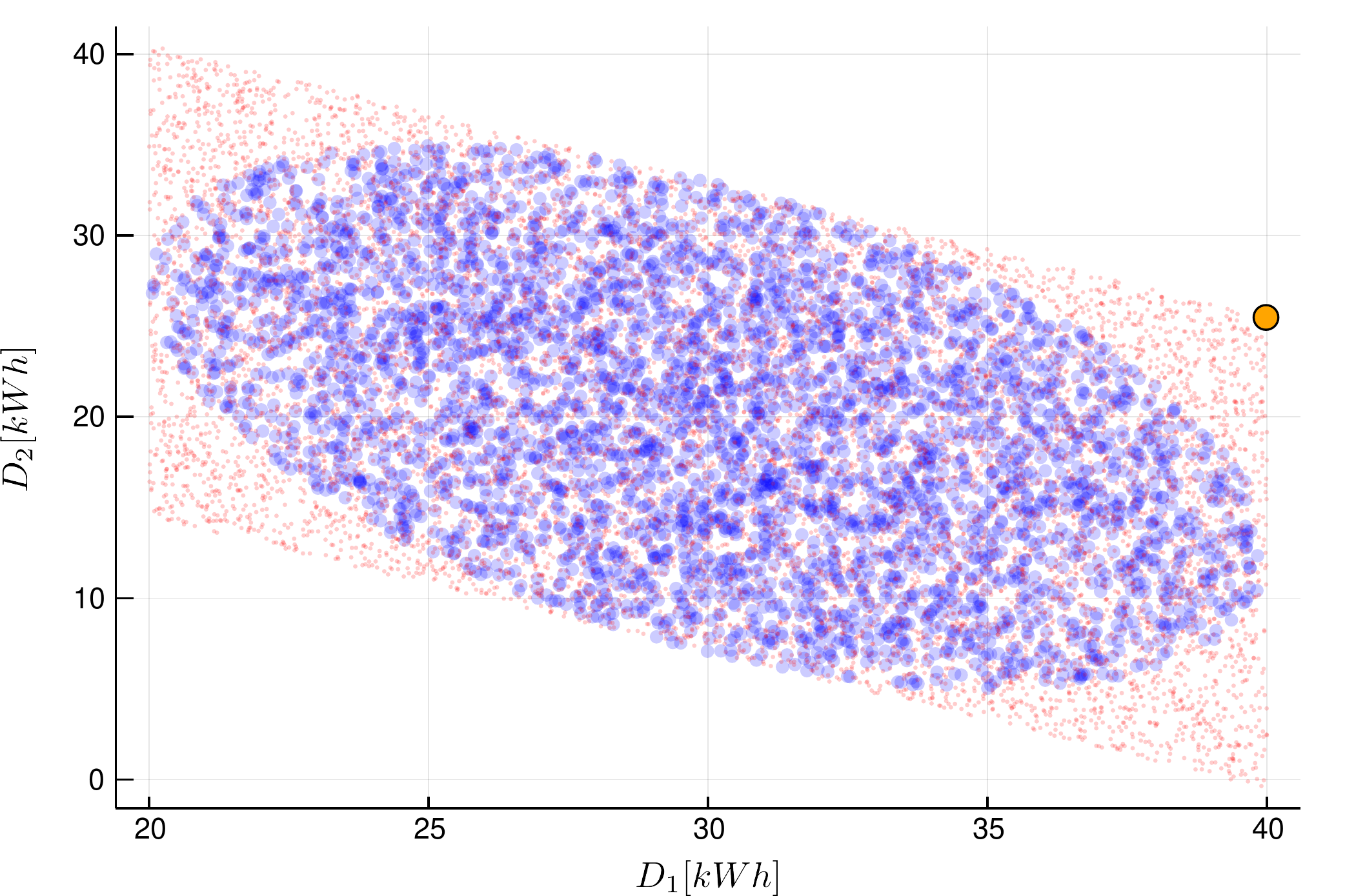}
\caption{Optimal total demand allocation (orange dot), joint bivariate distribution uncertainty set (blue dots), and Cholesky-based uncertainty approach (red dots) for correlation $-0.5$}
\label{fig:3}
\end{figure}

Figure \ref{fig:3} shows results for the case of $\rho = -0.5$. The optimal total allocation demand is represented by an orange dot. Uncertainty set is represented by blue dots for original ellipsoidal definition, and by red dots for the Cholesky-based linear approach used in our formulation.

\begingroup
\setlength{\tabcolsep}{10pt} 
\renewcommand{\arraystretch}{1.2} 
\begin{table}[t]
\centering
\caption{Results values for the total energy allocation, $D_i$, and the time-of-use of the loads, $L_i$, for different values of $\rho$.}

\begin{tabular}{d{2.1}|d{2.2} d{2.2}|d{2.2} d{2.2}|d{3.1}}
\hline
\multicolumn{1}{c|}{$\rho$}  & \multicolumn{1}{c}{$D_1$} & \multicolumn{1}{c|}{$D_2$} & \multicolumn{1}{c}{$L_1$} &  \multicolumn{1}{c|}{{$L_2$}} & \multicolumn{1}{c}{$C_T$ (\$)}  \\ \hline
-1 & 30.00 & 20.00 & 7.94 & 5.72 & 2610.5 \\
-0.5 & 40.00  & 25.48 & 11.40   & 9.60 & 3958.1 \\
0 & 40.00 & 35.00 & 10.55 & 10.00 & 4370.9 \\
+0.5 & 40.00 & 32.97 & 8.89  &  12.67 & 4159.8 \\
+1 & 30.00 &  20.00 & 7.94 &  5.72 & 2610.5 \\
\hline
\end{tabular}
\label{table:2}
\end{table}
\endgroup

\vspace{1cm}
\section{Conclusions}
We presented an alternative formulation for modeling time-of-use-type constraints of uninterruptible services. The proposed approach considers a single binary variable and the time of use parameter on the right-hand side of the allocation constraint. The formulation is useful for problems in which the time of use is uncertain. 
We applied the proposed model to the demand response problem for the allocation of loads. Results were illustrated in the context of deterministic and robust optimization framework.

\addtolength{\textheight}{-12cm}   






\bibliographystyle{IEEEtran}
\bibliography{references}


\end{document}